\documentclass[reqno,a4paper]{amsart}
\usepackage{amssymb}
\usepackage[colorlinks=true,hypertex]{hyperref}
\allowdisplaybreaks[4]
\theoremstyle{plain}
\newtheorem{thm}{Theorem}
\newtheorem{lem}{Lemma}
\newtheorem{open}{Problem}
\theoremstyle{remark}
\newtheorem{rem}{Remark}
\DeclareMathOperator{\td}{d\mspace{-2mu}}
\date{Completed on 5 September 2008, revised on 19 December 2008, in Melbourne}
\date{}

\begin{document}

\title[A double inequality for the arithmetic-geometric mean]
{An alternative and united proof of a double inequality for bounding the arithmetic-geometric mean}

\author[F. Qi]{Feng Qi}
\address[F. Qi]{Research Institute of Mathematical Inequality Theory, Henan Polytechnic University, Jiaozuo City, Henan Province, 454010, China}
\email{\href{mailto: F. Qi <qifeng618@gmail.com>}{qifeng618@gmail.com}, \href{mailto: F. Qi <qifeng618@hotmail.com>}{qifeng618@hotmail.com}, \href{mailto: F. Qi <qifeng618@qq.com>}{qifeng618@qq.com}}
\urladdr{\url{http://qifeng618.spaces.live.com}}

\author[A. Sofo]{Anthony Sofo}
\address[A. Sofo]{School of Engineering and Science, Victoria University, P.~O. Box 14428, Melbourne City, Victoria 8001, Australia}
\email{\href{mailto: A. Sofo <anthony.sofo@vu.edu.au>}{anthony.sofo@vu.edu.au}}

\begin{abstract}
In the paper, we provide an alternative and united proof of a double inequality for bounding the arithmetic-geometric mean.
\end{abstract}



\thanks{The first author was partially supported by the China Scholarship Council}

\thanks{This paper was typeset using \AmS-\LaTeX}

\maketitle

\section{Introduction}

The complete elliptic integral of the first kind was defined as
\begin{equation}
K(t)=\int_0^{\pi/2}\frac{\td\theta}{\sqrt{1-t^2\sin^2\theta}\,}
\end{equation}
for $0<t<1$, see~\cite[p.~132, Definition~3.2.1]{andrews}. It can also be defined in the following way: For positive numbers $a$ and $b$,
\begin{equation} \label{elliptic-1st-dfn-ab}
K(a,b)=\int_0^{\pi/2}\frac1{\sqrt{a^2\cos^2\theta+b^2\sin^2\theta}}\,\td \theta.
\end{equation}
\par
For positive numbers $a=a_0$ and $b=b_0$, let
\begin{equation}
a_{k+1}=\frac{a_k+b_k}2\quad \text{and} \quad b_{k+1}=\sqrt{a_kb_k}\,.
\end{equation}
In~\cite[p.~134, Definition~3.2.2]{andrews} and~\cite{Bracken-Expo-01}, the common limit $M(a,b)$ of these two sequences $\{a_k\}_{k=0}^\infty$ and $\{b_k\}_{k=0}^\infty$ is called as the arithmetic-geometric mean. It was proved in~\cite[Theorem~3.2.3]{andrews} and~\cite[Theorem~1]{Bracken-Expo-01} that
\begin{equation}
\frac1{M(a,b)}=\frac2\pi K(a,b).
\end{equation}
\par
For more information on the arithmetic-geometric mean and the complete elliptic integral of the first kind, please refer to~\cite[pp.~132--136]{andrews}, \cite{Bracken-Expo-01} and related references therein.
\par
In~\cite[Theorem~4]{Bracken-Expo-01} and~\cite{Carlson-Vuorinen-91}, it was proved that the inequality
\begin{equation}\label{M-L-ineq}
M(a,b)\ge L(a,b)
\end{equation}
holds true for positive numbers $a$ and $b$ and that the inequality~\eqref{M-L-ineq} becomes equality if and only if $a=b$, where
\begin{equation}
L(a,b)=\frac{b-a}{\ln b-\ln a}
\end{equation}
stands for the logarithmic mean for positive numbers $a$ and $b$ with $a\ne b$.
\par
In~\cite[Theorem~1.3]{Vamanamurthy-Vuorinen-JMAA-94}, it was turned out that
\begin{equation}\label{M-I-ineq}
M(a,b)<I(a,b)
\end{equation}
for positive numbers $a$ and $b$ with $a\ne b$, where
\begin{equation}
I(a,b)=\frac1e\biggl(\frac{b^b}{a^a}\biggr)^{1/(b-a)}
\end{equation}
represents the exponential mean for for positive numbers $a$ and $b$ with $a\ne b$.
\par
It is known that a generalization of the logarithmic mean $L(a,b)$ is the generalized logarithmic mean $L(p;a,b)$ of order $p\in\mathbb{R}$, which may be defined~\cite[p.~385]{bullenmean} by
\begin{equation}  \label{L(p;a,b)}
L(p;a,b)=
\begin{cases}
\biggl[\dfrac{b^{p+1}-a^{p+1}}{(p+1)(b-a)}\biggr]^{1/p}, & p\ne-1,0 \\
L(a,b), & p=-1 \\
I(a,b), & p=0
\end{cases}
\end{equation}
for positive numbers $a$ and $b$ with $a\ne b$, and that $L(p;a,b)$ is strictly increasing with respect to $p\in\mathbb{R}$. Therefore, one may naturally pose the following problem.

\begin{open}
What are the best constants $0\ge\beta>\alpha\ge-1$ such that the double inequality
\begin{equation}  \label{9.9-ineq-bracken-double}
L(\alpha;a,b)< M(a,b)< L(\beta;a,b)
\end{equation}
holds for all positive numbers $a$ and $b$ with $a\ne b$? In other words, are the constants $\alpha=-1$ and $\beta=0$ the best possible in the inequality~\eqref{9.9-ineq-bracken-double}?
\end{open}

It is easy to see that the complete elliptic integral $K(a,b)$ of the first kind tends to infinity as the ratio $\frac{b}{a}$ for $a>b>0$ tends to zero, equivalently, the arithmetic-geometric mean $M(a,b)$ tends to zero as $\frac{b}{a}\rightarrow0^{+}$. As $\frac{b}{a}\rightarrow0^{+}$, the logarithmic mean $L(a,b)$ also tends to $0$. However, the exponential mean $I(a,b)$ does not tend to zero as $\frac{b}{a}\rightarrow 0^{+}$. These phenomenons motivate us to put forward an alternative problem as follows.

\begin{open}\label{open-problem-ellip}
What are the best constants $\beta>\alpha\ge1$ such that the double inequality
\begin{equation}  \label{9.9-ineq-refine}
\alpha L(a,b)<M(a,b)<\beta L(a,b)
\end{equation}
holds for all positive numbers $a$ and $b$ with $a\ne b$?
\end{open}

In~\cite{sandor-jmaa-96} and~\cite[Theorem~1.3]{Vamanamurthy-Vuorinen-JMAA-94}, among others, the right-hand side inequality in~\eqref{9.9-ineq-refine} was verified to be valid for $\beta=\frac\pi2$.
\par
The aim of this paper is to confirm and sharpen the inequality~\eqref{9.9-ineq-refine} alternatively and unitedly.
\par
Our main result may be recited as the following theorem.

\begin{thm}\label{ellip-sharp-thm-2}
The double inequality~\eqref{9.9-ineq-refine} is valid for all positive numbers $a$ and $b$ with $a\ne b$ if and only if $\alpha\le1$ and $\beta\ge\frac\pi2$.
\end{thm}

\begin{rem}
Some inequalities were established in~\cite{Baricz-MZ-07, first-elliptic-integral.tex, ellip-xue-bao, construct, ellip-math-practice, ellip-gong-ke, ellip-huang, Yu-Kuang-Ye} for bounding the complete elliptic integrals.
\end{rem}

\section{Lemmas}

For proving our theorem alternatively and unitedly, we will employ the following lemmas.

\begin{lem}[\cite{Ponnusamy-Vuorinen-97}]\label{Ponnusamy-Vuorinen-97-lem}
Let $a_k$ and $b_k$ for $k\in\mathbb{N}$ be real numbers and the power series
\begin{equation}
A(x)=\sum_{k=1}^\infty a_kx^k\quad\text{and}\quad B(x)=\sum_{k=1}^\infty b_kx^k
\end{equation}
be convergent on $(-R,R)$ for some $R>0$. If $b_k>0$ and the ratio $\frac{a_k}{b_k}$ is strictly increasing for $k\in\mathbb{N}$, then the function $\frac{A(x)}{B(x)}$ is also strictly increasing on $(0,R)$.
\end{lem}

\begin{lem}[\cite{WallisFormula.html}]\label{wallis-formula}
For $n\in\mathbb{N}$,
\begin{equation}  \label{wsc}
\int_0^{\pi/2}\sin^nx\td x =\frac{\sqrt\pi\,\Gamma((n+1)/2)}{n\Gamma({n}/2)}
=\begin{cases}
\dfrac{\pi}2\cdot\dfrac{(n-1)!!}{n!!} & \text{for $n$ even}, \\[1em]
\dfrac{(n-1)!!}{n!!} & \text{for $n$ odd},
\end{cases}
\end{equation}
where $n!!$ denotes a double factorial.
\end{lem}

\begin{lem}\label{2-identities-sums}
For $k\in\mathbb{N}$,
\begin{gather}
\sum_{i=1}^{k}\binom{2i-2}{i-1}\frac{1}{i4^{i}} =\frac{1}{2}-\frac{2}{4^{k+1}}\binom{2k}{k},\label{h(k)}\\
\sum_{i=1}^{k}\binom{2i-2}{i-1}\frac{1}{(k-i+1)4^{i}} =\biggl[2\ln2+\gamma+\psi\biggl(k+\frac12\biggr)\biggr] \frac{\Gamma(k+1/2)}{4\sqrt{\pi}\,\Gamma(k+1)}.\label{g(k)}
\end{gather}
\end{lem}

\begin{proof}
For our own convenience, let us denote the two sequences in left-hand sides of~\eqref{h(k)} and~\eqref{g(k)} by $h(k)$ and $g(k)$ for $k\in\mathbb{N}$ respectively.
\par
When $k=1$, the identity~\eqref{h(k)} is valid clearly. Suppose the identity~\eqref{h(k)} holds for some $k>1$, then it follows that
\begin{align*}
h(k+1)&=h(k)+\binom{2k}{k}\frac1{(k+1)4^{k+1}}\\
&=\frac{1}{2}-\frac{2}{4^{k+1}}\binom{2k}{k} +\binom{2k}{k}\frac1{(k+1)4^{k+1}}\\
&=\frac{1}{2}-\frac{2}{4^{k+2}}\binom{2k+2}{k+1}.
\end{align*}
Therefore, by induction, the identity~\eqref{h(k)} is valid for all $k\in\mathbb{N}$.
\par
Applying the Zeilberger algorithm (see~\cite[Chapter~6]{a=b}) and~\eqref{wsc} yields
\begin{equation}\label{zeil-rec}
2(k+1)g(k+1)-(2k+1)g(k)=\frac{\Gamma(k+1/2)}{2\sqrt{\pi}\,\Gamma(k+1)} =\frac12\binom{2k}{k}\frac1{4^k}
\end{equation}
for $k\in\mathbb{N}$, from which the identity~\eqref{g(k)} follows.
\end{proof}

\begin{rem}
The identities~\eqref{h(k)} and~\eqref{g(k)} can also be verified easily by the famous softwares Maple and M\textsc{athematica}.
\end{rem}

\section{An alternative and united proof of Theorem~\ref{ellip-sharp-thm-2}}

Now we are in a position to alternatively and unitedly verify Theorem~\ref{ellip-sharp-thm-2}.
\par
Making use of the power series expansion
\begin{equation*}
\frac1{\sqrt{1-s}}=\sum_{i=0}^\infty\frac{(2i)!}{4^{i}(i!)^2}s^i,\quad 0<s<1,
\end{equation*}
it is obtained that
\begin{equation}  \label{series-ellip}
\frac1{\sqrt{1-s^2\sin^2\theta}}=\sum_{i=0}^\infty\frac{(2i)!}{4^{i}(i!)^2}s^{2i}\sin^{2i}\theta, \quad 0<s<1.
\end{equation}
From the celebrated Wallis sine formula~\eqref{wsc} in Lemma~\ref{wallis-formula}, it is obtained that
\begin{equation}  \label{2i-i-identity}
\int_0^{\pi/2}\sin^{2i}\theta\td \theta=\frac1{4^i}\binom{2i}{i}\frac\pi2 =\frac{\sqrt\pi\,}2\cdot\frac{\Gamma(i+1/2)}{\Gamma(i+1)},\quad i\in\mathbb{N}.
\end{equation}
Integrating on both sides of~\eqref{series-ellip} with respect to $\theta$ from $0$ to $\frac\pi2$ and using the identity~\eqref{2i-i-identity} yield
\begin{equation}  \label{series-fin}
\frac2\pi\int_0^{\pi/2}\frac{\td \theta}{\sqrt{1-s^2\sin^2\theta}}
=\sum_{i=0}^\infty\frac1{2^{4i}}{\binom{2i}{i}}^2s^{2i}
=\frac1\pi\sum_{i=0}^\infty\biggl[\frac{\Gamma(i+1/2)}{\Gamma(i+1)}\biggr]^2s^{2i}
\end{equation}
for $0<s<1$. Letting $s^2=1-t^2$ in~\eqref{series-fin} yields
\begin{equation}  \label{series-fin-t2}
\frac2\pi\int_0^{\pi/2}\frac{\td \theta}{\sqrt{1-(1-t^2)\sin^2\theta}}
=\sum_{i=0}^\infty\frac1{2^{4i}}{\binom{2i}{i}}^2\bigl(1-t^2\bigr)^{i}
\triangleq \sum_{k=0}^\infty b_k\bigl(1-t^2\bigr)^{k}
\end{equation}
for $0<t<1$.
\par
It is easy to see that
\begin{equation}
\biggl(\frac1t\biggr)^{(k)}=\frac{(-1)^kk!}{t^{k+1}}\quad \text{and}\quad
\bigl(\sqrt{t}\,\bigr)^{(k)}=\frac{(-1)^{k+1}(2k-1)!!}{2^k(2k-1)t^{k-1/2}}
\end{equation}
for $k\in\mathbb{N}$. Then, by Leibniz's Theorem~\cite[p.~12, 3.38]{abram} for differentiation of a product, we gain that
\begin{multline*}
\bigl[\bigl(1+\sqrt{t}\,\bigr)\ln t\bigr]^{(k)}=\bigl(1+\sqrt{t}\,\bigr)
\biggl(\frac1t\biggr)^{(k-1)} +\bigl(1+\sqrt{t}\,\bigr)^{(k)}\ln t \\
+\sum_{i=1}^{k-1}\binom{k}{i}\bigl(1+\sqrt{t}\,\bigr)^{(i)}\biggl(\frac1t\biggr)^{(k-i-1)} =\bigl(1+\sqrt{t}\,\bigr)\frac{(-1)^{k-1}(k-1)!}{t^{k}} \\
+\bigl(\sqrt{t}\,\bigr)^{(k)}\ln t +\frac{(-1)^kk!}{t^{k-1/2}}
\sum_{i=1}^{k-1}\frac{(2i-1)!!}{(2i)!!(2i-1)(k-i)},
\end{multline*}
where, and elsewhere in this paper, an empty sum is understood to be nil. Thus,
\begin{equation}
\bigl[\bigl(1+\sqrt{t}\,\bigr)\ln t\bigr]^{(k)}\Big\vert_{t=1} =(-1)^kk!%
\Biggl[\sum_{i=1}^{k-1}\frac{(2i-1)!!}{(2i)!!(2i-1)(k-i)}-\frac2k\Biggr]
\end{equation}
for $k\in\mathbb{N}$. Hence,
\begin{equation}
\bigl(1+\sqrt{t}\,\bigr)\ln t=\sum_{k=1}^\infty \Biggl[\sum_{i=1}^{k-1}\frac{(2i-1)!!}{(2i)!!(2i-1)(k-i)}-\frac2k\Biggr](1-t)^k
\end{equation}
which can be reduced by replacing $\sqrt{t}\,$ by $t$ to
\begin{equation}
(1+t)\ln t=\sum_{k=1}^\infty \Biggl[\frac12\sum_{i=1}^{k-1} \frac{(2i-1)!!}{(2i)!!(2i-1)(k-i)}-\frac1k\Biggr]\bigl(1-t^2\bigr)^k,
\end{equation}
and so
\begin{multline*}
\frac{\ln t}{t-1}=\frac{(1+t)\ln t}{t^2-1} =\sum_{k=1}^\infty\Biggl[\frac1k-\frac12\sum_{i=1}^{k-1} \frac{(2i-1)!!}{(2i)!!(2i-1)(k-i)}\Biggr]\bigl(1-t^2\bigr)^{k-1} \\
=\sum_{k=0}^\infty \Biggl[\frac1{k+1}-\frac12\sum_{i=1}^{k} \frac{(2i-1)!!}{(2i)!!(2i-1)(k-i+1)}\Biggr] \bigl(1-t^2\bigr)^{k} \triangleq\sum_{k=0}^\infty a_k\bigl(1-t^2\bigr)^{k}
\end{multline*}
for $0<t<1$.
\par
The two identities in Lemma~\ref{2-identities-sums} and the equality in the right-hand side of the inequality~\eqref{2i-i-identity} give
\begin{align*}
\frac{1}{k+1}-a_{k}& =\frac{1}{2}\sum_{i=1}^{k}\frac{(2i-1)!}{2^{i-1}(i-1)!2^{i}i!(2i-1)(k-i+1)} \\
& =\sum_{i=1}^{k}\binom{2i-2}{i-1}\frac{1}{2^{2i}i(k-i+1)} \\
& =\frac{1}{k+1}\sum_{i=1}^{k}\binom{2i-2}{i-1}\frac{1}{4^{i}}\biggl(\frac{1}{i}+\frac{1}{k-i+1}\biggr) \\
& =\frac{1}{k+1}\sum_{i=1}^{k}\binom{2i-2}{i-1}\frac{1}{4^{i}i}+\frac{1}{k+1} \sum_{i=1}^{k}\binom{2i-2}{i-1}\frac{1}{4^{i}(k-i+1)} \\
& =\frac{1}{k+1}\biggl[\frac{1}{2}-\frac{2}{4^{k+1}}\binom{2k}{k}\biggr] +\frac{1}{k+1}\sum_{i=1}^{k}\binom{2i-2}{i-1}\frac{1}{4^{i}(k-i+1)} \\
& =\frac{1}{k+1}\biggl\{\frac{1}{2}-\frac{2}{4^{k+1}}\binom{2k}{k}+\frac{[2\ln 2 +\gamma+\psi(k+1/2)]\Gamma(k+1/2)}{4\sqrt{\pi }\,\Gamma (k+1)}\biggr\} \\
& =\frac{1}{k+1}\biggl[\frac{1}{2}-\frac{2-2\ln 2-\gamma -\psi (k+1/2)}{4\sqrt{\pi }\,}\cdot \frac{\Gamma(k+1/2)}{\Gamma (k+1)}\biggr],
\end{align*}
that is,
\begin{equation*}
a_{k}=\frac{1}{k+1}\biggl[\frac{1}{2}+\frac{2-2\ln 2-\gamma-\psi (k+1/2)}{4\sqrt{\pi }\,}\cdot \frac{\Gamma(k+1/2)}{\Gamma (k+1)}\biggr],
\end{equation*}
where $k\in \mathbb{N}$ and $\gamma =0.57721566\dotsm$ is Euler-Mascheroni's constant. It is listed in~\cite[p.~258, 6.3.4]{abram} that
\begin{equation}
\psi \biggl(n+\frac{1}{2}\biggr)=-\gamma -2\ln 2+2\biggl(1+\frac{1}{3}
+\dotsm +\frac{1}{2n-1}\biggr),\quad n\geq 1.
\end{equation}%
Hence,
\begin{equation}
a_{k}=\frac{1}{2(k+1)}\Biggl[1-\frac{\Gamma (k+1/2)}{\sqrt{\pi}\,\Gamma(k+1)} \Biggl(\sum_{i=1}^{k}\frac{1}{2i-1}-1\Biggr)\Biggr],\quad k\in \mathbb{N}.
\end{equation}
\par
Now let us discuss the increasingly monotonic property of the ratio $\frac{a_k}{b_k}$ for $k\in\mathbb{N}$. It is clear that
\begin{equation}  \label{ratio-a-k}
\frac{a_k}{b_k}\le\frac{a_{k+1}}{b_{k+1}}\quad \Longleftrightarrow \quad
\frac{b_{k+1}}{b_k}\le\frac{a_{k+1}}{a_k} \quad \Longleftrightarrow \quad
\biggl(\frac{2k+1}{2k+2}\biggr)^2\le\frac{a_{k+1}}{a_k}
\end{equation}
which is equivalent to
\begin{equation}  \label{rearra-19}
\Biggl[\sum_{i=2}^{k+1}\frac1{2i-1}-\frac{(k+1/2)(k+2)}{(k+1)^2}
\sum_{i=2}^{k}\frac1{2i-1}\Biggr] \frac{\Gamma(k+3/2)}{\sqrt\pi\,\Gamma(k+2)}
\le1-\frac{(k+1/2)^2(k+2)}{(k+1)^3}
\end{equation}
for $k\ge2$. Furthermore, easy calculation gives
\begin{align*}
&\quad\sum_{i=2}^{k+1}\frac1{2i-1}-\frac{(k+1/2)(k+2)}{(k+1)^2}\sum_{i=2}^{k}\frac1{2i-1} \\
&=\frac1{2k+1}+\biggl[1-\frac{(k+1/2)(k+2)}{(k+1)^2}\biggr]\sum_{i=2}^{k}\frac1{2i-1} \\
&=\frac1{2k+1}-\frac{k}{2(k+1)^2}\sum_{i=2}^{k}\frac1{2i-1} \\
&=\frac{k}{2(k+1)^2}\Biggl[\frac{2(k+1)^2}{k(2k+1)}-\sum_{i=2}^{k}\frac1{2i-1}\Biggr].
\end{align*}
Since the sequence $\frac{2(k+1)^2}{k(2k+1)}$ for $k\ge2$ is strictly decreasing and tends to $1$ as $k\to\infty$ and the sequence $\sum_{i=2}^{k}\frac1{2i-1}$ for $k\ge2$ is strictly increasing and diverges to $\infty$, the sequence
\begin{equation}  \label{seq-frac-sum}
S_k\triangleq\frac{2(k+1)^2}{k(2k+1)}-\sum_{i=2}^{k}\frac1{2i-1}
\end{equation}
for $k\ge2$ is strictly decreasing and diverges to $-\infty$. As a result, there exists an integer $k_0\ge2$ such that the sequence $S_k$ is negative for all $k\ge k_0$. From the fact that $S_9=0.01\dotsm$ and $S_{10}=-0.04\dotsm$, it follows that $k_0=10$. Therefore, considering the facts that
\begin{equation*}
\frac{(k+1/2)^2(k+2)}{(k+1)^3}\le1
\end{equation*}
and
\begin{equation*}
\frac{\Gamma(k+3/2)}{\sqrt\pi\,\Gamma(k+2)}>0
\end{equation*}
for $k\ge2$, it readily follows that the inequality~\eqref{rearra-19} holds for all $k\ge10$.
\par
Straightforward computations reveal that
\begin{center}
\begin{tabular}{|c||c|c|c|c|c|c|c|c|c|c|}
\hline
$k$ & $1 $ & $2 $ & $3 $ & $4 $ & $5 $ & $6 $ & $7 $ & $8 $ & $9 $ & $10 $
\\ \hline
$a_k$ & $\frac{1}{4} $ & $\frac{7}{48} $ & $\frac{5}{48} $ & $\frac{313}{3840} $ & $\frac{43}{640} $ & $\frac{12317}{215040} $ & $\frac{10751}{215040}$ & $\frac{183349}{4128768} $ & $\frac{206329}{5160960} $ & $\frac{66087019}{1816657920} $ \\ \hline
\end{tabular}
\end{center}
and that
\begin{center}
\begin{tabular}{|c||c|c|c|}
\hline
$k$ & $1$ & $2$ & $3$ \\ \hline
$\frac{a_{k+1}}{a_k}$ & $\frac{7}{12}=0.583\dotsm$ & $\frac{5}{7}=0.714\dotsm
$ & $\frac{313}{400}=0.782\dotsm$ \\ \hline
$\bigl(\frac{2k+1}{2k+2}\bigr)^2$ & $\frac{9}{16}=0.562\dotsm$ & $\frac{25}{%
36}=0.694\dotsm$ & $\frac{49}{64}=0.765\dotsm$ \\ \hline\hline
$k$ & $4$ & $5$ & $6$ \\ \hline
$\frac{a_{k+1}}{a_k}$ & $\frac{258}{313}=0.824\dotsm$ & $\frac{12317}{14448}%
=0.852\dotsm$ & $\frac{10751}{12317}=0.872\dotsm$ \\ \hline
$\bigl(\frac{2k+1}{2k+2}\bigr)^2$ & $\frac{81}{100}=0.810\dotsm$ & $\frac{121%
}{144}=0.840\dotsm$ & $\frac{169}{196}=0.862\dotsm$ \\ \hline\hline
$k$ & $7$ & $8$ & $9$ \\ \hline
$\frac{a_{k+1}}{a_k}$ & $\frac{916745}{1032096}=0.888\dotsm$ & $\frac{825316%
}{916745}=0.900\dotsm$ & $\frac{66087019}{72627808}=0.909\dotsm$ \\ \hline
$\bigl(\frac{2k+1}{2k+2}\bigr)^2$ & $\frac{225}{256}=0.878\dotsm$ & $\frac{%
289}{324}=0.892\dotsm$ & $\frac{361}{400}=0.902\dotsm$ \\ \hline
\end{tabular}
\end{center}
Consequently, the inequality~\eqref{rearra-19} holds for all $k\in\mathbb{N}$. This demonstrates that the ratio $\frac{a_k}{b_k}$ for $k\in\mathbb{N}$ is strictly increasing. By virtue of Lemma~\ref{Ponnusamy-Vuorinen-97-lem}, it is obtained that the function
\begin{equation}  \label{ratio-means-two-t}
\frac{(\ln t)/(t-1)}{(2/\pi)\int_0^{\pi/2}1\Big/\sqrt{1-(1-t^2)\sin^2\theta}%
\,\td\theta}
\end{equation}
is strictly decreasing in $t\in(0,1)$. The well-known H\^ospital's rule yields that the limits of the function~\eqref{ratio-means-two-t} are $\frac\pi2$ and $1$ as $t$ tends to $0^+$ and $1^-$ respectively. Hence, the double inequality
\begin{equation}  \label{double-ineq-2-mean}
\frac2\pi\cdot\frac{\ln t}{t-1}<\frac2\pi\int_0^{\pi/2}\frac{\td\theta}{\sqrt{1-(1-t^2)\sin^2\theta}\,} <\frac{\ln t}{t-1}
\end{equation}
for $t\in(0,1)$ is valid and sharp. Letting $t=\frac{a}b$ for $b>a>0$ in~\eqref{double-ineq-2-mean} leads to
\begin{equation}  \label{double-ineq-2-m-ab}
\frac2\pi\cdot\frac{\ln a-\ln b}{a-b} <\frac2\pi\int_0^{\pi/2} \frac{\td \theta} {\sqrt{a^2\sin^2\theta+b^2\cos^2\theta}\,} <\frac{\ln a-\ln b}{a-b}.
\end{equation}
It is easy to see that the inequality~\eqref{double-ineq-2-m-ab} is valid for all positive numbers $a$ and $b$ with $a\ne b$. This implies that the double inequality~\eqref{9.9-ineq-refine} is valid for all positive numbers $a$ and $b$ with $a\ne b$ if and only if $\alpha\le\frac2\pi$ and $\beta\ge1$. The proof of Theorem~\ref{ellip-sharp-thm-2} is complete.

\subsection*{Acknowledgements}
The authors appreciate Associate Professor Chao-Ping Chen in China for his recommending the reference~\cite{Bracken-Expo-01} and some valuable comments in an e-mail to the first author on 8 June 2008. This paper was completed while the first author was visiting the Research Group in Mathematical Inequalities and Applications, Victoria University, Melbourne, Australia between March 2008 and February 2009 by the grant from the China Scholarship Council. So the first author expresses many thanks to the local colleagues such as Professors Pietro Cerone and Sever S. Dragomir for their invitation and hospitality throughout this period.

\end{document}